\documentclass[a4paper,11pt]{article}
\usepackage{amssymb}
\usepackage{amsfonts}
\usepackage{amsmath}

\input{tcilatex}

\begin{document}

\title{Topology of unitary groups and the prime orders of binomial
coefficients\\
{\small Dedicated to Professor Boju Jiang on his 80th birthday}}
\author{Haibao Duan\thanks{%
Supported by NSFC 11131008; 11661131004} \\
Institute of Mathematics, Chinese Academy of Sciences,\\
School of Mathematical Sciences, \\
University of the Chinese Academy of Sciences \and Xianzu Lin\thanks{%
Supported by National Natural Science Foundation for young (no.11401098).} \\
College of Mathematics and Computer Science, \\
Fujian Normal University, Fuzhou, 350108, China;}
\maketitle

\begin{abstract}
Let $c:SU(n)\rightarrow PSU(n)=SU(n)/\mathbb{Z}_{n}$ be the quotient map of
the special unitary group $SU(n)$ by its center subgroup $\mathbb{Z}_{n}$.
We determine the induced homomorphism $c^{\ast }:$ $H^{\ast
}(PSU(n))\rightarrow H^{\ast }(SU(n))$ on cohomologies by computing with the
prime orders of binomial coefficients.

\begin{description}
\item[2000 Mathematical Subject Classification: ] 55T10

\item[Key words and phrases:] Lie groups; cohomology, prime order of an
integer

\item[Email:] dhb@math.ac.cn;linxianzu@126.com
\end{description}
\end{abstract}

\section{Main result}

The center subgroup of the special unitary group $SU(n)$ is the cyclic group 
$\mathbb{Z}_{n}$ of order $n$ generated by the diagonal matrix $diag\{e^{%
\frac{2\pi }{n}i},\cdots ,e^{\frac{2\pi }{n}i}\}\in SU(n)$. Let $%
c:SU(n)\rightarrow PSU(n):=SU(n)/\mathbb{Z}_{n}$ be the quotient
homomorphism. The quotient group $PSU(n)$ is called \textsl{the projective
unitary group} of rank $n$. Consider the induced map of $c$ on cohomologies

\begin{quote}
$c^{\ast }:H^{\ast }(PSU(n);R)\rightarrow H^{\ast }(SU(n);R)$ with $R=$ $%
\mathbb{Q}$ or $\mathbb{Z}$.
\end{quote}

\noindent It is well known that if $R=\mathbb{Q}$ then $c^{\ast }$ is an
isomorphism of algebras. It is also known that if $R=\mathbb{Z}$ there are
integral classes $\xi _{2r-1}\in H^{2r-1}(SU(n);\mathbb{Z})$, $2\leq r\leq n$%
, so that the integral cohomology $H^{\ast }(SU(n);\mathbb{Z})$ is the
exterior ring $\Lambda (\xi _{3},\cdots ,\xi _{2n-1})$. It follows that
there exist integral cohomology classes

\begin{quote}
$\zeta _{2r-1}\in H^{2r-1}(PSU(n);\mathbb{Z})$, $2\leq r\leq n$,
\end{quote}

\noindent as well as a sequence $(a_{n,2},\cdots ,a_{n,n})$ of integers,
such that

\begin{enumerate}
\item[(1.1)] the set $\{1,\zeta _{I}=\underset{r\in I}{\Pi }\zeta
_{2r-1}\mid $ $I\subseteq \{2,\cdots ,n\}\}$ of square free monomials in the 
$\zeta _{2r-1}$'s is a basis of the free part of the ring $H^{\ast }(PSU(n);%
\mathbb{Z})$;

\item[(1.2)] the induced map of the covering $c$ satisfies the relation

$c^{\ast }(\zeta _{2r-1})=a_{n,r}\cdot \xi _{2r-1}$, $2\leq r\leq n$.
\end{enumerate}

\noindent In addition, since the mapping degree of $c$ is the order of the
center $\mathbb{Z}_{n}$, and since the product $\zeta _{3}\cdots \zeta
_{2n-1}$ is a generator of the top degree cohomology $H^{m}(PSU(n);\mathbb{Z}%
)=\mathbb{Z}$ ($m=\dim PSU(n)$) by (1.1), the calculation by (1.2)

\begin{quote}
$c^{\ast }(\zeta _{3}\cdots \zeta _{2n-1})=(a_{n,2}\cdots a_{n,n})\cdot \xi
_{3}\cdots \xi _{2n-1}$
\end{quote}

\noindent indicates that the sequence $(a_{n,2},\cdots ,a_{n,n})$ of
integers must satisfy the constraint

\begin{enumerate}
\item[(1.3)] $n=a_{n,2}\cdot \cdots \cdot a_{n,n}$.
\end{enumerate}

\noindent In view of (1.2) and (1.3) it is both of topological and
arithmetic interests to study the problem of expressing the sequence $%
(a_{n,2},\cdots ,a_{n,n})$ as explicit function in $n$.

By \textsl{the prime factorization }of an integer $n\geq 2$ we mean the
unique expression $n=p_{1}^{r_{1}}\cdots p_{t}^{r_{t}}$, where $%
1<p_{1}<\cdots <p_{t}$ is the set of all prime factors of $n$. In term of
this factorization we introduce the partition on the set $\left\{ 2,\cdots
,n\right\} $ by

\begin{quote}
$\left\{ 2,\cdots ,n\right\} =Q_{0}(n)\underset{1\leq i\leq t}{\amalg }%
Q_{p_{i}}(n)$,
\end{quote}

\noindent where $Q_{p_{i}}(n)=\{p_{i}^{r}\mid 1\leq r\leq r_{i}\}$, and
where $Q_{0}(n)$ is the complement of the subset $\underset{1\leq i\leq t}{%
\amalg }Q_{p_{i}}(n)\subset \left\{ 2,\cdots ,n\right\} $. The main result
of this paper is

\bigskip

\noindent \textbf{Theorem 1.1. }\textsl{One has}\textbf{\ }$a_{n,k}=p_{i}%
\QTR{sl}{\ }$\textsl{or }$1$ \textsl{in accordance to}\textbf{\ }$k\in
Q_{p_{i}}(n)$ $\QTR{sl}{or}$ $k\in Q_{0}(n)$\textsl{.}

\bigskip

The plan of this paper is: based on an arithmetic characterization of the
integers $a_{n,k}$ obtained in Lemma 2.2, we establish Theorem 1.1 in
Section 3. In this paper the cohomologies are over the ring $\mathbb{Z}$ of
integers, unless otherwise stated.

The authors would like to thank their referees for useful comments and
suggestions.

\section{Preliminaries in topology}

The center of the unitary group $U(n)$ of order $n$ is the circle subgroup $%
S^{1}=\{diag\{e^{i\theta },\cdots ,e^{i\theta }\}\mid $ $\theta \in \lbrack
0,2\pi ]\}$ that acts on $U(n)$ via right multiplication. The quotient group 
$U(n)/S^{1}$ can be identified with the projective unitary group $PSU(n)$ 
\cite{BB}, while the quotient map $C:U(n)\rightarrow PSU(n)$ is both a group
homomorphism and an oriented circle bundle over $PSU(n)$ that is related to
the covering $c$ by the commutative diagram

\begin{center}
\begin{tabular}{lll}
&  & $U(n)$ \\ 
& $\overset{i}{\nearrow }$ & $\quad \downarrow C$ \\ 
$SU(n)$ & $\overset{c}{\rightarrow }$ & $PSU(n)$%
\end{tabular}%
,
\end{center}

\noindent where $i$ is the obvious inclusion. Let $\omega \in H^{2}(PSU(n))$
be the Euler class of the circle bundle $C:U(n)\rightarrow PSU(n)$. Then the
Gysin sequence \cite[p.149]{MS} provides us with an exact sequence relating
the cohomologies $H^{\ast }(U(n))$ and $H^{\ast }(PSU(n))$

\begin{enumerate}
\item[(2.1)] ${\small \cdots \rightarrow H}^{r}{\small (PSU(n))}\overset{%
{\small C}^{\ast }}{{\small \rightarrow }}{\small H}^{r}{\small (U(n))}%
\overset{\theta }{{\small \rightarrow }}{\small H}^{r-1}{\small (PSU(n))}%
\overset{{\small \cup \omega }}{{\small \rightarrow }}{\small H}^{r+1}%
{\small (PSU(n))}\overset{{\small C}^{\ast }}{{\small \rightarrow }}{\small %
\cdots }$.
\end{enumerate}

\noindent It has been shown in \cite[Formula (4.10)]{D} that, with respect
to the presentation $H^{\ast }(U(n))=$\ $\Lambda (\xi _{1},\cdots ,\xi
_{2n-1})$\ by Borel \cite{B}, one has

\bigskip

\noindent \textbf{Lemma 2.1.} \textsl{The connecting homomorphism }$\theta $%
\textsl{\ in (2.1) satisfies the relation}

\begin{enumerate}
\item[(2.2)] $\theta (\xi _{2k-1})=\binom{n}{k}\omega ^{k-1}$.$\square $
\end{enumerate}

For a pair $(n,k)$ of integer with $1\leq k\leq n$ we set

\begin{quote}
$b_{n,k}:=g.c.d.\{\binom{n}{1},\cdots ,\binom{n}{k}\}$,
\end{quote}

\noindent where $\binom{n}{k}=\frac{n!}{k!(n-k)!}$ is the binomial
coefficient. Note that $b_{n,k}\mid b_{n,k-1}$ with $b_{n,1}=n$ and $%
b_{n,n}=1$. The following result gives an arithmetic characterization of the
integers $a_{n,k}$ in (1.2).

\bigskip

\noindent \textbf{Lemma 2.2. }\textsl{There exists a set }$\Phi =\{\zeta
_{2k-1}\in H^{2k-1}(PSU(n))\mid $\textsl{\ }$2\leq k\leq n\}$\textsl{\ of
cohomology classes satisfying the following properties}

\begin{enumerate}
\item[(2.3)] $c^{\ast }(\zeta _{2k-1})=\frac{b_{n,k-1}}{b_{n,k}}\cdot \xi
_{2k-1}$, $2\leq k\leq n$;

\item[(2.4)] \textsl{the set }$\Psi =\{1,\zeta _{I}=\underset{r\in I}{\Pi }%
\zeta _{2r-1}\mid $\textsl{\ }$I\subseteq \{2,\cdots ,n\}\}$\textsl{\ of
square free monomials is a basis of the free part of the group }$H^{\ast
}(PSU(n))$\textsl{.}
\end{enumerate}

\noindent \textsl{In particular, one has} $a_{n,k}=\frac{b_{n,k-1}}{b_{n,k}}$%
\textsl{.}

\bigskip

\noindent \textbf{Proof. }The determinantal function $\det :U(n)\rightarrow
S^{1}$ is both a group homomorphism and a smooth submersion

\begin{quote}
$SU(n)\overset{i}{\hookrightarrow }U(n)\overset{\det }{\rightarrow }S^{1}$
(i.e. with fiber the subgroup $SU(n)$).
\end{quote}

\noindent It gives rise to a retraction $p:$ $U(n)\rightarrow SU(n)$ by

\begin{quote}
$p(g)=diag\{\det (g)^{-1},1,\cdots ,1\}\cdot g$, $g\in U(n)$.
\end{quote}

\noindent It follows from $i^{\ast }\circ p^{\ast }=id:$ $H^{\ast
}(SU(n))\rightarrow $ $H^{\ast }(SU(n))$ that the map $p^{\ast }$ carries $%
H^{\ast }(SU(n))$ isomorphically onto the first summand of the obvious
decomposition

\begin{enumerate}
\item[(2.5)] $H^{\ast }(U(n))=\Lambda (\xi _{3},\cdots ,\xi _{2n-1})\oplus
\xi _{1}\cdot \Lambda (\xi _{3},\cdots ,\xi _{2n-1})$,
\end{enumerate}

\noindent while $i^{\ast }$ maps the first summand isomorphically onto $%
H^{\ast }(SU(n))$.

By the exactness of the Gysin sequence (2.1) one gets from (2.2) that

\begin{quote}
$\omega \cup \theta (\xi _{2k-1})=\binom{n}{k}\omega ^{k}=0$ for all $1\leq
k\leq n$.
\end{quote}

\noindent Consequently,

\begin{quote}
$b_{n,k-1}\cdot \theta (\xi _{2k-1})=0$ for all $2\leq k\leq n$.
\end{quote}

\noindent In particular, for any $k>1$ the order of the class $\theta (\xi
_{2k-1})\in H^{\ast }(PSU(n))$ divides $\frac{b_{n,k-1}}{b_{n,k}}$. By the
exactness of the sequence (2.1) there exists a set $\Phi =\{\zeta _{2k-1}\in
H^{2k-1}(PSU(n))\mid $ $2\leq k\leq n\}$ of cohomology classes satisfying

\begin{quote}
$C^{\ast }(\zeta _{2k-1})=\frac{b_{n,k-1}}{b_{n,k}}\cdot \xi _{2k-1}$, $%
2\leq k\leq n$.
\end{quote}

\noindent Applying $i^{\ast }$ to both sides one obtains formula (2.3) by
the relation $i^{\ast }\circ C^{\ast }=c^{\ast }$ and by (2.5).

Setting $m=\dim PSU(n)$ the monomial $\xi _{3}\cdots \xi _{2n-1}$ is a
generator of the top degree cohomology group $H^{m}(SU(n))=\mathbb{Z}$ by
(2.5). Since the mapping degree of $c$ is $n$ the calculation

\begin{quote}
$c^{\ast }(\zeta _{3}\cdots \zeta _{2n-1})=b_{n,1}\cdot (\xi _{3}\cdots \xi
_{2n-1})=n\cdot (\xi _{3}\cdots \xi _{2n-1})$
\end{quote}

\noindent by (2.3) indicates that the monomial $\zeta _{3}\cdots \zeta
_{2n-1}$ is a generator of the top degree cohomology group $H^{m}(PSU(n))=%
\mathbb{Z}$. As a result the set $\Psi $ of $2^{n-1}$ monomials in (2.4) is
linearly independent in $H^{\ast }(PSU(n))$.

We claim further that the elements in $\Psi $ spans a direct summand of the
group $H^{\ast }(PSU(n))$. Assume, on the contrary, that there exist a
monomial $\zeta _{I}\in \Psi $, a class $\varsigma \in H^{\ast }(PSU(n))$,
as well as some integer $a>1$, so that a relation of the form $\zeta
_{I}=a\cdot \varsigma $ holds in $H^{\ast }(PSU(n))$. Multiplying both sides
by the class $\zeta _{\overline{I}}$ with $\overline{I}$ the complement of $%
I\subseteq \{2,\cdots ,n\}$ yields the equality

\begin{quote}
$\zeta _{3}\cdots \zeta _{2n-1}=(-1)^{r}a\cdot (\varsigma \cup \zeta _{%
\overline{I}})$ (for some $r\in \mathbb{Z}$).
\end{quote}

\noindent This contradicts to the fact that $\zeta _{3}\cdots \zeta _{2n-1}$
generates $H^{m}(PSU(n))=\mathbb{Z}$. The proof of (2.4) is completed by the
routine relation

\begin{quote}
$\dim (H^{\ast }(PSU(n))\otimes \mathbb{Q)=}2^{n-1}$.$\square $
\end{quote}

\section{The prime orders of binomial coefficients}

Let $\mathbb{R}^{+}$ and $\mathbb{N}$ be, respectively, the set of positive
reals and the set of natural numbers. For a real $x\in \mathbb{R}$ let $[x]$
denote the unique integer satisfying $0\leq x-[x]<1$. It is straightforward
to see that

\bigskip

\noindent \textbf{Lemma 3.1.} \textsl{For any} $x,y\in \mathbb{R}^{+}$, $%
m,n\in \mathbb{N}$, \textsl{with} $x+y=m+n$,\textsl{\ we have }

\begin{quote}
\textsl{\qquad }$[x]+[y]\leq \lbrack m]+[n],$
\end{quote}

\noindent \textsl{where the inequality holds} \textsl{if} \textsl{either} $%
x\neq \lbrack x]$\textsl{\ or }$y\neq \lbrack y]$\textsl{.}$\square $

\bigskip

Given a prime $p$ and an $m\in \mathbb{N}$ \textsl{the order of }$m$ \textsl{%
at} $p$, denoted by $ord_{p}(m)$, is the biggest integer $a$ so that $m$ is
divisible by the power $p^{a}$. Clearly, if $n_{1},\cdots ,n_{r}$ is a
sequence of positive integers then

\begin{enumerate}
\item[(3.1)] $ord_{p}(g.c.d.\{n_{1},\cdots
,n_{r}\})=min\{ord_{p}(n_{1}),\cdots ,ord_{p}(n_{r})\}.$
\end{enumerate}

\noindent The next formula is shown in \cite[Theorem 416]{HW}

\bigskip

\noindent \textbf{Lemma 3.2.} $ord_{p}(n!)=\sum_{k=1}^{\infty }[\tfrac{n}{%
p^{k}}]$.

\bigskip

Applying Lemma 3.2 we show that

\bigskip

\noindent \textbf{Lemma 3.3.} \textsl{For an }$k\in Q_{p}(n)$\textsl{\ the
number }$a_{n,k}=\frac{b_{n,k-1}}{b_{n,k}}$\textsl{\ is divisible by }$p$%
\textsl{.}

\bigskip

\noindent \textbf{Proof. }In view of the formula (3.1) it suffices to show,
for any $1\leq k<p^{s}$ and $1\leq s\leq r$, that

\begin{quote}
$ord_{p}(\binom{n}{p^{s}})<ord_{p}(\binom{n}{k})$.
\end{quote}

\noindent This is equivalent to

\begin{enumerate}
\item[(3.2)] $ord_{p}(p^{s}!(n-p^{s})!)>ord_{p}(k!(n-k)!)$, $1\leq k<p^{s}$, 
$1\leq s\leq r$.
\end{enumerate}

\noindent By Lemma 3.2, if $1\leq k\leq p^{s}$, then

\begin{enumerate}
\item[(3.3)] $ord_{p}(k!(n-k)!)$ $=\sum_{j=1}^{\infty }([\tfrac{n-k}{p^{j}}%
]+[\tfrac{k}{p^{j}}])=\sum_{j=1}^{s}([\tfrac{n-k}{p^{j}}]+[\tfrac{k}{p^{j}}%
])+\sum_{j=s+1}^{\infty }[\tfrac{n-k}{p^{j}}]$,
\end{enumerate}

\noindent where the second equality comes from the obvious relation

\begin{quote}
$\sum_{j=s+1}^{\infty }[\tfrac{k}{p^{j}}]=0$.
\end{quote}

\noindent Since the integer $n$ is divisible by $p^{s}$ the open interval $(%
\tfrac{n-p^{s}}{p^{j}},\tfrac{n}{p^{j}})$ with $j>s$ contains no integer
(for otherwise, there would be an integer $m\in (0,p^{s})$ such that $n-m$
is divisible by the power $p^{j}$ for some $j>s$, contradicting to that $m$
is not divisible by $p^{s}$). It follows that

\begin{quote}
$[\tfrac{n-p^{s}}{p^{j}}]=[\tfrac{n-k}{p^{j}}]$ for all $1\leq k<p^{s}$, $%
j>s $.
\end{quote}

\noindent Summing this equality over $j\geq s+1$ yields that

\begin{enumerate}
\item[(3.4)] $\sum_{j=s+1}^{\infty }[\tfrac{n-k}{p^{j}}]=\sum_{j=s+1}^{%
\infty }[\tfrac{n-p^{s}}{p^{j}}].$
\end{enumerate}

\noindent Finally, the relation (3.2) is shown by the calculation

\begin{quote}
$ord_{p}(p^{s}!(n-p^{s})!)=\sum_{j=1}^{s}([\tfrac{n-p^{s}}{p^{j}}]+[\tfrac{%
p^{s}}{p^{j}}])+\sum_{j=t+1}^{\infty }[\tfrac{n-p^{s}}{p^{j}}]$ (by (3.3))

$\quad >\sum_{j=1}^{s}([\tfrac{n-k}{p^{j}}]+[\tfrac{k}{p^{j}}%
])+\sum_{j=s+1}^{\infty }[\tfrac{n-p^{s}}{p^{j}}]$ (by Lemma 3.1)

$\quad =\sum_{j=1}^{s}([\tfrac{n-k}{p^{j}}]+[\tfrac{k}{p^{j}}%
])+\sum_{j=s+1}^{\infty }[\tfrac{n-k}{p^{j}}]$ (by (3.4))

$\quad =ord_{p}(k!(n-k)!)$ (by (3.3)).$\square $
\end{quote}

We come now to a proof of Theorem 1.1 stated in Section 1.

\bigskip

\noindent \textbf{Proof of Theorem 1.1. }Assume that the prime factorization
of the integer $n\geq 2$ is $n=p_{1}^{r_{1}}\cdots p_{t}^{r_{t}}$. For each $%
k\in Q_{p_{i}}(n)$ we can assume by Lemma 3.3 that $a_{n,k}=p_{i}\cdot 
\overline{a}_{n,k}$ for some $\overline{a}_{n,k}\in \mathbb{N}$, and
rephrase the decomposition (1.3) of the integer $n$ as

\begin{quote}
$n=p_{1}^{r_{1}}\cdots p_{t}^{r_{t}}\underset{k\in Q_{p_{i}}(n),1\leq i\leq t%
}{\Pi }\overline{a}_{n,k}\underset{k\in Q_{0}(n)}{\Pi }a_{n,k}$.
\end{quote}

\noindent Since the factor $p_{1}^{r_{1}}\cdots p_{t}^{r_{t}}$ on the right
hand side coincides with the prime factorization of $n$, the equality above
forces out the relations $\overline{a}_{n,k}=1$ for $k\in Q_{p_{i}}(n)$, and 
$a_{n,k}=1$ for $k\in Q_{0}(n)$. This completes the proof.$\square $

\bigskip

Let $\Delta (\zeta _{3},\cdots ,\zeta _{2n-1})\subset H^{\ast }(PSU(n))$ be
the subgroup spanned additively by the elements in $\Psi $, and let $T$ the
torsion ideal of the ring $H^{\ast }(PSU(n))$. By Lemma 2.2 the cohomology $%
H^{\ast }(PSU(n))$ admits the decomposition

\begin{enumerate}
\item[(3.5)] $H^{\ast }(PSU(n))=\Delta (\zeta _{3},\cdots ,\zeta
_{2n-1})\oplus T$.
\end{enumerate}

\noindent Since $H^{\ast }(SU(n))$ is torsion free $c^{\ast }(T)=0$. Assume
that the prime factorization of the integer $n>2$ is $p_{1}^{r_{1}}\cdots
p_{t}^{r_{t}}$. Theorem 1.1 implies, with respect to the decomposition
(3.5), the following result:

\bigskip

\noindent \textbf{Theorem 3.4.} \textsl{The induced map }$c^{\ast }:$\textsl{%
\ }$H^{\ast }(PSU(n))\rightarrow H^{\ast }(SU(n))$\textsl{\ is}

\begin{quote}
$c^{\ast }(T)=0$ \textsl{and} $c^{\ast }(\zeta _{I})=p_{1}^{s_{1}}\cdots
p_{t}^{s_{t}}\cdot \xi _{I}$\textsl{,}
\end{quote}

\noindent \textsl{where} $I\subseteq \{2,\cdots ,n\}$, $s_{i}=\left\vert
I\cap Q_{p_{i}}(n)\right\vert $, $1\leq i\leq t$.$\square $

\bigskip

For another application of Theorem 1.1 consider a compact Lie group $G$ with
a maximal torus $T$. Let $\pi ^{\ast }:H^{\ast }(G/T)\rightarrow H^{\ast
}(G) $ be the induced map of the torus fibration $\pi :G\rightarrow G/T$. By
Grothendieck \cite{G} the subring $\func{Im}\pi ^{\ast }$ $\subset H^{\ast
}(G)$ is the \textsl{Chow ring} $A^{\ast }(G^{c})$ of the reductive
algebraic group $G^{c}$ corresponding to $G$.

Let $J_{n}(\omega )\subset H^{\ast }(PSU(n))$ be the subring generated by $%
\omega $. The proofs of Theorem 1.1 and Lemma 2.2 indicates that

\begin{enumerate}
\item[(3.6)] $J_{n}(\omega )=\mathbb{Z}[\omega ]/\left\langle b_{n,k}\omega
^{k},1\leq k\leq n\right\rangle $, $\deg \omega =2$.
\end{enumerate}

\noindent On the other hand one can show for $G=PSU(n)$ that

\begin{enumerate}
\item[(3.7)] $\func{Im}\pi ^{\ast }=J_{n}(\omega )$.
\end{enumerate}

\noindent Moreover, in view of the obvious decomposition $b_{n,k}=a_{n,k+1}$ 
$\cdots a_{n,n}$ and by Theorem 1.1, the number $b_{n,k}$ admits the prime
factorization

\begin{enumerate}
\item[(3.8)] $b_{n,k}=\underset{1\leq i\leq t}{\Pi }p_{i}^{r_{i}-s_{i}}$,%
\textsl{\ }
\end{enumerate}

\noindent where $s_{i}=r_{i}$ if $k\geq p_{i}^{r_{i}}$, and satisfies the
relation $s_{i}\leq \log _{p_{i}}k<s_{i}+1$\ if $k<p_{i}^{r_{i}}$. Combining
(3.6), (3.7) with (3.8) yields the following decomposition of the Chow ring $%
A^{\ast }(PSU(n)^{c})$ into its primary torsion ideals

\bigskip

\noindent \textbf{Theorem 3.5.} \textsl{If }$n>2$ \textsl{has} \textsl{the
prime factorization }$p_{1}^{r_{1}}\cdots p_{t}^{r_{t}}$\textsl{, then }

\begin{quote}
$A^{\ast }(PSU(n)^{c})=\mathbb{Z}\underset{i\in \{1,\cdots ,t\}}{\mathbb{%
\oplus }}$ $\mathbb{Z}[\omega ]^{+}/\left\langle p_{i}^{r_{i}}\omega
,p_{i}^{r_{i}-1}\omega ^{p_{i}},p_{i}^{r_{i}-2}\omega ^{p_{i}^{2}},\cdots
,\omega ^{p_{i}^{r_{i}}}\right\rangle $.$\square $
\end{quote}

\noindent \textbf{Remark 1.} In \cite[\S 4]{BB} Baum and Browder obtained
certain information on the Serre spectral sequence associated to the fiber
sequence $U(n)\overset{C}{\rightarrow }PU(n)\overset{\chi }{\rightarrow }%
BS^{1}$, where $\chi $ is the classifying map of the Euler class $\omega \in
H^{2}(PU(n))$. Write $H^{\ast }(U(n))=\Lambda (z_{i}:i=1,\cdots ,n)$ and $%
H^{\ast }(BS^{1})=\mathbb{Z}[\alpha ]$, where $\left\vert z_{i}\right\vert
=2i-1$ and $\left\vert \alpha \right\vert =2$, It is shown in \cite[\S 4]{BB}
that the differentials $d_{r}$ in the spectral sequence satisfy the relations

\begin{enumerate}
\item[(3.8)] $d_{r}(z_{i})$ $=$ $0$ for $r<2i$; $d_{2i}(z_{i})$ $=$ $\binom{n%
}{i}\alpha ^{i}$.
\end{enumerate}

\noindent Applying (3.8) one can obtain an alternative proof of Lemma 2.2.

Our approach to Lemmas 2.2 uses the exact sequence (2.1) which can be
extended to compute the cohomologies of all the adjoint Lie groups, e.g. 
\cite[Theorem 4.7, Theorem 4.12]{D}.

\bigskip

\noindent \textbf{Remark 2}. The study of the prime orders of binomial
coefficients has a long and outstanding history, see the survey article \cite%
{Gra} by Granville. Using the classical Kummer's Theorem \cite{Gra} one can
give an alternative proof of Lemma 3.3. We note that the statement and
expression of Kummer's result is too technical to be presented here, while
the our proof of Lemma 3.3 is independent to Kummer's approach, and is more
accessible to the general audience.

\end{document}